\documentclass[a4paper,10pt]{article}
\usepackage[all]{xy}
\usepackage{graphicx}
\usepackage[T1]{fontenc}
\usepackage[utf8]{inputenc}
\usepackage{amsmath,amsthm}
\usepackage{amsfonts,amssymb}
\usepackage{bigdelim}
\usepackage{multirow}
\usepackage{enumerate}
\usepackage{float}
\usepackage{pstricks}
\usepackage[ae]{}
\usepackage{authblk}
\usepackage[numbers]{natbib}
\newcommand{\E}{\mathbb{E}}
\newcommand{\N}{\mathbb{N}}

\newcommand{\R}{\mathbb{R}}

\newcommand{\Pb}{\mathbb{P}}

\def\={{\;\mathop{=}\limits^{\text{(law)}}\;}}

\newtheorem{theorem}{Theorem}[section]

\newtheorem{defi}[theorem]{Definition}

\theoremstyle{definition}
\newtheorem{rem}[theorem]{Remark}

\newtheorem{exa}[theorem]{Example}

\numberwithin{equation}{section}

\usepackage{fancyhdr}
\date{}

\usepackage[english]{babel}

\title{Self-similar martingales derived from Root embedding}
\author{\Large Antoine-Marie Bogso\footnote{University of Yaounde I,
 Department of Mathematics, P.O. Box 812 Yaoundé, Cameroon. Phone: (+237)652 620 452, Email: ambogso@gmail.com} and Mohamed Mbehou\footnote{University of Yaounde I,
 Department of Mathematics, P.O. Box 812 Yaoundé, Cameroon. Phone: (+237)675518266, Email: mbehoumoh@gmail.com } }
\date{}

\begin{document}
 
\maketitle
%\faketableofcontents
\begin{abstract}Given a family $(\mu_\lambda,\lambda\geq0)$ of integrable mean-zero probability measures such that, for every $\lambda\geq0$, $\mu_\lambda$ is the image of $\mu_1$ under the homothety $y\longmapsto\sqrt{\lambda}y$, we provide a necessary and sufficient condition on $\mu_1$ under which the Root embedding algorithm yields a self-similar martingale with one-dimensional marginals $(\mu_\lambda,\lambda\geq0)$. Precisely, if $\tau_{\lambda}$ and $R_{\lambda}$ denote the Root solution to the Skorokhod embedding problem (SEP) and the Root regular barrier for $\mu_\lambda$ respectively, then this condition is equivalent to the property that $(R_{\lambda},\lambda\geq0)$ is non-increasing in the sense of inclusion, which in turn is equivalent to the assertion that $(\tau_\lambda,\lambda\geq0)$ is non-decreasing a.s. We show that there are many examples for which this result applies and we provide some numerical simulations to illustrate the monotonicity property of regular barriers $(R_{\lambda},\lambda\geq0)$ in this case. 
\text{}
  \\
\textbf{keywords: } Skorokhod embedding problem, Root embedding, regular barriers, self-similar martingales.\\
 \textbf{subclass MSC:} 60E15, 60G44, 60J25.
\end{abstract}

\section{Introduction}
\label{intro}
There are many results in the literature related to the construction of self-similar processes. In Madan-Yor \cite{MY02}, Fan-Hamza-Klebaner \cite{FHK15}, Hamza-Klebaner \cite{HK07}, Hirsch-Profeta-Roynette-Yor \cite{HPRY11}, Bogso \cite{Bo15} and Henry-Labordère-Tan-Touzi \cite{HLTT16}, the authors provided many constructions of self-similar martingales with given marginal distributions. In particular, Madan and Yor \cite{MY02}, Hamza and Klebaner \cite{HK07}, and Henry-Labord\`ere, Tan and Touzi \cite{HLTT16} exhibited several examples of discontinuous fake Brownian motions. Albin \cite{Al08} answered positively the question of the existence of continuous fake Brownian motion, and this result was extended by Baker-Donati-Martin-Yor \cite{BDY11} who exhibited a sequence of continuous martingales with Brownian marginal distributions and scaling property. A quite simple construction of continuous fake Brownian motion, based on Box-Muller transform, has been given by Oleszkiewicz \cite{Ol08}. The results on fake Brownian motion was extended by Hobson \cite{Ho13} to prove the existence of a continuous fake exponential Brownian motion. More recently, Jourdain and Zhou \cite{JZ16} provided a new class of fake Brownian motions which solve a special class of local and stochastic volatility  SDEs. Certain of the works cited above use Skorokhod embedding solutions to construct self-similar martingales. Precisely, Madan and Yor \cite{MY02} exploit the Az\'ema-Yor algorithm, Hirsch, Profeta, Roynette and Yor \cite{HPRY11} apply Az\'ema-Yor, Hall-Breiman and Bertoin-Le Jan embedding solutions, and they provided a new Skorokhod embedding solution that gave another class of self-similar martingales. We show that the Root embedding solution also provides a class of martingales which enjoy Brownian scaling.

Let $\mu$ be a square-integrable mean-zero probability measure and let $(B_t,t\geq0)$ be a standard Brownian motion. Root \cite{Ro69} proved the existence of a closed time-space set $R$, the so-called {\it Root barrier, }such that the first hitting time $\tau$ of $R$ by the time-space process $(t,B_t;t\geq0)$ solves the Skorokhod embedding problem for $\mu$ (SEP($\mu$)), meaning that $B_\tau$ has distribution $\mu$ and $(B_{\tau\wedge t},t\geq0)$ is uniformly integrable. He also defined a barrier function $r:\,[-\infty,+\infty]\to[0,+\infty]$ attached to $R$ as $r(x)=\inf\{t\geq0:\,(t,x)\in R\}$ and observed that $r$ is lower semi-continuous. The problem of the existence of a stopping-time for Brownian motion in such a way that the stopped value has a given distribution was first stated and solved by Skorokhod \cite{Sk14}. Note that different  Root barriers may embed the same distribution. This was solved by Loynes \cite{Lo70} who introduced the notion of a {\it regular }barrier and proved that there exists exactly one regular barrier that solves SEP($\mu$). On the other hand, 
the Root's solution is optimal in the sense that it has minimal variance among all stopping times $S$ such that $B_S$ has distribution $\mu$ and $\E[S]=\int y^2\mu(dy)$. This was conjectured by Kiefer \cite{Kie72} and was solved later by  Rost \cite{Rst76}. We refer to Beiglb{\"o}ck, Cox and 
Huesmann \cite{BCH17} where a transport-based approach to the SEP has been developed to derive all known and a variety of new optimal solutions.
Another interesting question on Root embedding is that it is not easy to find explicitely the regular Root barrier for a given distribution. Indeed, this barrier is constructed explicitely only for a handful of simple examples. Dupire \citep{Du05} showed formally that the Root barrier $R$ is given by the solution of a nonlinear PDE. This was further developed by Cox and Wang \cite{CW13} who use a variational formulation to calculate $R$. A complete characterization of regular Root barriers as free boundaries of PDEs has been provided by Gassiat, Mijatovic and Dos Reis \cite{GODR17}. When the distribution $\mu$ is atom-free, Gassiat, Mijatovic and Oberhauser \cite{GMO15} established  that the barrier function $r$ solves a nonlinear Volterra integral equation and that if, in addition, $r$ is continuous, then $r$ is the unique solution. More recently, Cox, Obl\`oj and Touzi \cite{COT17} provided a characterization of regular Root barriers by means of an optimal stopping formulation and exploited this approach to establish a finitely-many marginals extension of the Root solution to the SEP. These authors also proved that their solution satisfies an optimality property which extends the optimality property of the one-marginal Root solution. Using the Cox, Hobson and Touzi results, Richard, Tan and Touzi \cite{RTT18} provided a full marginals extension on some compact time interval of the Root solution to SEP. Precisely, using a tightness result established by K{\"a}llblad, Tan and Touzi \cite[Lemma 4.5]{KTT17}, they proved that the full marginals limit of the finitely-many marginals Root solution for the SEP exists and enjoys the same optimality property as the multiple-marginals Root solution provided in \cite{COT17}. 

Here we consider the case of a family $(\mu_\lambda,\lambda\geq0)$ of integrable mean-zero probability measures where $\mu_\lambda$ is the image of $\mu_1$ under the homothety $y\longmapsto\sqrt{\lambda}y$. 
We apply Root solution for the SEP to embed simultaneously all $\mu_\lambda$'s into a standard Brownian motion $(B_t,t\geq0)$ issued from $0$. The Root embedding provides a family of regular barriers $(R_\lambda,\lambda\geq0)$ and a family of stopping times $(\tau_\lambda,\lambda\geq0)$ such that $\tau_\lambda=\inf\{t\geq0;\,(t,B_t)\in R_\lambda\}$ and $B_{\tau_\lambda}$ has law $\mu_\lambda$. We apply the optimal stopping characterization of one-marginal Root solution to the SEP given in \cite[Theorem 2.8]{COT17} and Brownian scaling to prove that the regular barrier function $(\lambda,x)\longmapsto r_\lambda(x):=\inf\{t\geq0:\,(t,x)\in R_{\lambda}\}$ defined on $[0,+\infty]\times[-\infty,+\infty]$ is self-similar in the sense that 
$$
r_{\lambda}(x)=\lambda r_1\left(\dfrac{x}{\sqrt{\lambda}}\right),\quad\forall\,(\lambda,x)\in\R_+^*\times[-\infty,+\infty],
$$ 
where $\R_+^*$ denotes the set of positive real numbers. This result can also be deduced from the viscosity PDE characterization of regular Root barriers obtained by Gassiat, Mijatovic and Dos Reis \cite[Theorem 2]{GODR17}.
We deduce from the monotonicity property of $r$ that, for every $x\in\R$, $\lambda\longmapsto r_\lambda(x)$ is non-decreasing. The self-similarity property of the function $(\lambda,x)\longmapsto r_{\lambda}(x)$ given above allows us to obtain a necessary and sufficient condition on $r_1$ under which the family of Root barriers $(R_\lambda,\lambda\geq0)$ is non-increasing, in the sense that $R_\lambda \subset R_\delta$ (i.e. $R_\delta$ includes $R_\lambda$) for every $0\leq\delta\leq\lambda$. This monotonicity property of the family $(R_\lambda,\lambda\geq0)$ is equivalent to the assertion that $\lambda\longmapsto\tau_\lambda$ is non-decreasing a.s. Then, as $\tau_\lambda$ solves the SEP for $\mu_\lambda$, $\left(B_{\tau_\lambda},\lambda\geq0\right)$ is a martingale with marginals $(\mu_\lambda,\lambda\geq0)$. Moreover, we prove that $\left(B_{\tau_\lambda},\lambda\geq0\right)$ enjoys the Brownian scaling and the Markovian properties.

In Section 2, we exploit the optimal stopping characterization of one-marginal Root solution for the SEP to prove that the function $(\lambda,x)\longmapsto r_\lambda(x)$ is self-similar. Then, we provide a sufficient condition on the barrier function $r_1$ under which the family of regular Root barrier $(R_\lambda,\lambda\geq0)$ is non-increasing. This allows us to exhibit a new class of martingales with Brownian scaling in Section 3. We also discuss the Markovian properties of these processes. In Section 4, we provide some numerical simulations to illustrate the monotonicity property of the regular barriers $(R_\lambda,\lambda\geq0)$. The numerical scheme follows the idea of Gassiat-Oberhauser-Dos Reis \cite[Section 4]{GODR17}. In particular, the Barles-Souganidis method \cite{BS90,BS91} can be applied to obtain  the convergence of the scheme and a result due to Jakobsen \cite{Ja03} provides its convergence rate.

\section{Root embedding under scaling}
Let $\mu$ be an integrable probability measure, and let $(B_t,t\geq0)$ denote a one-dimensional Brownian. A solution to SEP($\mu$) is any stopping time $\tau$ such that $B_{\tau}$ has law $\mu$, and $\left(B_{\tau\wedge t},t\geq0\right)$ is uniformly integrable. In the case where $\mu$ has zero mean and a second moment, Root provided a Skorokhod embedding solution  that is the first hitting time of a barrier, the so-called {\it Root barrier}.
\begin{defi}\label{defi:RootBarrier}A closed subset $R$ of $[0,+\infty]\times[-\infty,+\infty]$ is called a Root barrier if 
\begin{enumerate}
\item[(i)]$(t,x)\in R$ implies $(s,x)\in R$ for all $s\geq t$,
\item[(ii)]$(+\infty,x)\in R$ for all $x\in [-\infty,+\infty]$,
\item[(iii)]$[0,+\infty]\times\{-\infty,+\infty\}\subset R$.
\end{enumerate}
\end{defi}
Given  a Root barrier, one defines its barrier function $r:[-\infty,+\infty]\to[0,+\infty]$ as
$$
r(x)=\inf\{t\geq0:\,(t,x)\in R\},\qquad x\in[-\infty,+\infty].
$$
Since $R$ is closed, then, as observed by  Root \citep{Ro69} and Loynes \citep{Lo70}, $r$ is a lower semi-continuous function. Moreover, one deduces from Property (i) in Defition \ref{defi:RootBarrier} that $R$ is the epigraph of $r$ in the plane $[0,+\infty]\times[-\infty,+\infty]$ (see e.g. Cox-Obl\`oj-Touzi \citep{COT17}), i.e.
$$
R=\{(t,x)\in[0,+\infty]\times[-\infty,+\infty]:\,t\geq r(x)\}.
$$

There may exist different Root barriers which solve SEP($\mu$).   But Loynes \citep{Lo70} introduced the notion of {\it regular barrier} and he provided a uniqueness result when we restrict ourselves to regular Root barriers.  
\begin{defi}
A Root barrier $R$ is said to be regular if its barrier function $r$ vanishes outside the interval $[x_{-},x_{+}]$, where $x_{-}$ and $x_{+}$ are respectively the first negative and the first positive zeros of $r$.
\end{defi}

\begin{theorem}(Loynes \citep{Lo70}, Rost \citep{Rst76}, Gassiat-Mijatovic-Oberhauser \citep{GMO15}). Suppose that $\mu$ is of finite variance and has zero mean. Then there exists exactly one regular Root barrier  $R$ such that $\tau=\inf\{t\geq0:\,(t,B_t)\in R\}$ solves SEP($\mu$). Moreover, $\tau$ minimises for every $t\geq0$ the residual expectation $\E\left[(\widetilde{\tau}-t)^{+}\right]$ among all $\widetilde{\tau}$ that are solutions of SEP($\mu$).
\end{theorem}
The finite variance assumption in the preceding result has recently been relaxed to the condition that the measure has a finite first moment. This was first obtained by Gassiat-Oberhauser-Dos Reis \citep{GODR17} who provide a complete characterization of regular Root barriers as free boundaries of PDEs. The next result is a special case of Theorem 2 and Corollary 1 in \citep{GODR17}. 
\begin{theorem}\label{theo:RootFreeBounds}(Gassiat-Oberhauser-Dos Reis \citep[Theorem 2, and Corollary 1]{GODR17}). Let $\mu$ be an integrable and centered probability measure. For every $\lambda\geq0$, let $\mu_\lambda$ be the image measure of $\mu$ under $y\longmapsto\sqrt{\lambda}y$ (in particular, $\mu=\mu_1$). The following equivalent assertions hold:
\begin{enumerate}
\item[(i)]There exists a regular Root barrier $R_{\lambda}$ such that $\tau_{\lambda}=\inf\{t\geq0:\,(t,B_t)\in R_{\lambda}\}$ solves SEP($\mu_{\lambda}$),
\item[(ii)]There exists a viscosity solution $u_{\lambda}\in C([0,+\infty],[-\infty,+\infty])$, decreasing in time, of
\begin{equation}\label{eq:RootPDE}
\left\{
\begin{array}{ll}
\min\left(u-v_{\mu_{\lambda}}, \partial_tu-\dfrac{1}{2}\partial_{xx} u\right)=0&\text{on }]0,+\infty[\times\R,\\
u(0,\cdot)=-|\cdot|&\text{on }\R,\\
u(+\infty,\cdot)=v_{\mu_{\lambda}} &\text{on }\R,
\end{array}
\right.
\end{equation}
where $v_{\mu_{\lambda}}$ is the potential function of $\mu_{\lambda}$: $$v_{\mu_{\lambda}}(x)=-\int_{\R}|x-y|\mu_{\lambda}(dy),\quad x\in\R.$$ 
\end{enumerate}
Moreover,
\begin{equation}\label{eq:RootFreeBounds}
R_{\lambda}=\{(t,x)\in[0,+\infty]\times[-\infty,+\infty]:\,u_{\lambda}(t,x)=v_{\mu_{\lambda}}(x)\}\text{ and } u_{\lambda}(t,x)=-\E\left[\left|B_{\tau_{\lambda}\wedge t}-x\right|\right].
\end{equation}
Precisely, $R_{\lambda}$ is the unique regular Root barrier such that $\tau_{\lambda}$ solves SEP($\mu_\lambda$) and $R_{\lambda}$ is the free boundary (\ref{eq:RootFreeBounds}) of the obstacle PDE (\ref{eq:RootPDE}).
\end{theorem}
\begin{rem}\label{rem:RootBounds}
Gassiat-Oberhauser-Dos Reis \citep{GODR17} also provided a comparison theorem that allows them to prove that $u_\lambda$ is the unique viscosity solution of linear growth of (\ref{eq:RootPDE}) (see \citep[Theorem 5]{GODR17}).
\end{rem}
A characterization of Root solution to the Skorokhod embedding problem by means of an optimal stopping formulation has been provided recently by Cox-Obl\'oj-Touzi \cite[Theorem 2.8]{COT17}. These authors proved this result using purely probabilistic methods. We state here a special case of Theorem 2.8 in \citep{COT17}.
\begin{theorem}\label{theo:OptStopRoot}(Cox-Obl\'oj-Touzi \cite[Theorem 2.8]{COT17}). 
Consider a one-dimensional Brownian motion $(B_s,s\geq0)$ defined on a filtered probability space $(\Omega,\mathcal{F},(\mathcal{F}_s,s\geq0),\Pb)$ satisfying the usual conditions.  Let $\mu$ be an integrable zero-mean probability measure. For every $\lambda\geq0$, let $\mu_\lambda$ denote the image of $\mu$ under $y\longmapsto\sqrt{\lambda}y$. Define 
\begin{equation}\label{eq:RootOptStopSol}
u_\lambda(t,x)=\sup\limits_{\sigma\in\mathcal{T}_B^t}J^{\lambda,B}_{t,x}(\sigma)\quad\text{with}
\end{equation}
\begin{equation}\label{eq:OptJlambda}
J^{\lambda,B}_{t,x}(\sigma)=\E\left[v_{\delta_0}(x+B_\sigma)+\left(v_{\mu_\lambda}-v_{\delta_0}\right)(x+B_\sigma)\mathbf{1}_{\sigma<t}\right],
\end{equation}
where $\mathcal{T}_B^t$ is the collection of all $(\mathcal{F}_s,s\geq0)$-stopping times $\sigma$ taking values in $[0,t]$. Then the stopping region
\begin{equation}\label{eq:RootStopReg}
R_\lambda=\left\{(t,x)\in[0,+\infty]\times[-\infty,+\infty]:\,u_\lambda(t,x)=v_{\mu_\lambda}(x)\right\}
\end{equation}
is the regular barrier inducing the Root solution to the SEP($\mu_\lambda$). 
Moreover,
$$
u_\lambda(t,x)=-\E\left[\left\vert B_{\tau_\lambda\wedge t}-x\right\vert\right].
$$
\end{theorem} 
For every $\lambda>0$, let $r_\lambda$ denotes the barrier function of $R_\lambda$. The next result states that the function $(\lambda,x)\longmapsto r_\lambda(x)$ satisfies a self-similarity property.   
\begin{theorem}\label{theo:RootPDEscale}
Let $\mu$ be an integrable and centered probability measure. For every $\lambda\geq0$, let $\mu_\lambda$ denote the image of $\mu$ under $y\longmapsto\sqrt{\lambda}y$. Let $u_\lambda$ be the function defined by (\ref{eq:RootOptStopSol}) and (\ref{eq:OptJlambda}). Then,
\begin{enumerate}
\item[(i)]for every $\lambda>0$,
$$
u_\lambda(t,x)=\sqrt{\lambda}u_1\left(\frac{t}{\lambda},\frac{x}{\sqrt{\lambda}}\right),\quad\forall\,(t,x)\in[0,+\infty]\times[-\infty,+\infty],
$$
where $u_1$ is given by
\begin{equation*} 
u_1(t,x)=\sup\limits_{\sigma\in\mathcal{T}_B^t}J^{1,B}_{t,x}(\sigma)\quad\text{with}
\end{equation*}
\begin{equation*}
J^{1,B}_{t,x}(\sigma)=\E\left[v_{\delta_0}(x+B_\sigma)+\left(v_{\mu_1}-v_{\delta_0}\right)(x+B_\sigma)\mathbf{1}_{\sigma<t}\right],
\end{equation*} 
\item[(ii)]the map $(\lambda,x)\longmapsto r_\lambda(x)$ is self-similar in the sense that
\begin{equation}\label{eq:BarrierFunctEq}
r_\lambda(x)=\lambda r_1\left(\frac{x}{\sqrt{\lambda}}\right),\quad\forall\,(\lambda,x)\in\R_+^\ast\times[-\infty,+\infty].
\end{equation}
\end{enumerate}
\end{theorem}
 
\begin{proof} 
\item[(i)]By Brownian scaling, 
$$
\left(W_s=\frac{1}{\sqrt{\lambda}}B_{\lambda s},\mathcal{G}_s=\mathcal{F}_{\lambda s}; s\geq0\right)
$$
is a Brownian motion. As a consequence, $u_1$ rewrites
$$
u_1(t,x)=\sup\limits_{\rho\in\mathcal{T}_W^t}J^{1,W}_{t,x}(\rho)\quad\text{with}
$$
$$
J^{1,W}_{t,x}(\rho)=\E\left[v_{\delta_0}(x+W_\rho)+\left(v_{\mu_1}-v_{\delta_0}\right)(x+W_\rho)\mathbf{1}_{\rho<t}\right],
$$
where $\mathcal{T}_{W}^t$ is the collection of all $(\mathcal{G}_{s},s\geq0)$-stopping times $\rho$ which take values in $[0,t]$.   Moreover, for every $\lambda>0$, $(t,x)\in\R_+\times\R$, and $\sigma\in\mathcal{T}^t_B$,
\begin{align*}
J^{\lambda,B}_{t,x}(\sigma)&=\E\left[v_{\delta_0}(x+B_\sigma)+\left(v_{\mu_\lambda}-v_{\delta_0}\right)(x+B_\sigma)\mathbf{1}_{\sigma<t}\right]\\
&=\sqrt{\lambda}\,\E\left[v_{\delta_0}\left(\frac{x+B_\sigma}{\sqrt{\lambda}}\right)+\left(v_{\mu_1}-v_{\delta_0}\right)\left(\frac{x+B_\sigma}{\sqrt{\lambda}}\right) \mathbf{1}_{\sigma<t}\right]\\
&=\sqrt{\lambda}\,\E^{}\left[v_{\delta_0}\left(\frac{x}{\sqrt{\lambda}}+W_{\sigma/\lambda}\right)+\left(v_{\mu_1}-v_{\delta_0}\right)\left(\frac{x}{\sqrt{\lambda}}+W_{\sigma/\lambda}\right)\mathbf{1}_{\sigma/\lambda<t/\lambda}\right]\\&=\sqrt{\lambda}J^{1,W}_{\frac{t}{\lambda},\frac{x}{\sqrt{\lambda}}}(\sigma/\lambda).
\end{align*}
Hence, as $(1/\lambda)\mathcal{T}^t_B:=\left\{\sigma/\lambda;\,\sigma\in\mathcal{T}^t_B\right\}=\mathcal{T}^{\frac{t}{\lambda}}_W$, we have
$$
u_\lambda(t,x)=\sup\limits_{\sigma\in\mathcal{T}^t_B}J^{\lambda,B}_{t,x}(\sigma)=\sqrt{\lambda}\sup\limits_{\rho\in\mathcal{T}_W^{\frac{t}{\lambda}}}J^{1,W}_{\frac{t}{\lambda},\frac{x}{\sqrt{\lambda}}}(\rho)=\sqrt{\lambda}\,u_1\left(\frac{t}{\lambda},\frac{x}{\sqrt{\lambda}}\right).
$$
\item[(ii)] We deduce from Point (i) that
\begin{align*}
R_\lambda&=\{(t,x)\in[0,+\infty]\times[-\infty,+\infty]:\,u_\lambda(t,x)=v_{\mu_\lambda}(x)\}\\
&=\left\{(t,x)\in[0,+\infty]\times[-\infty,+\infty]:\,u_1\left(\frac{t}{\lambda},\frac{x}{\sqrt{\lambda}}\right)=v_{\mu_1}\left(\frac{x}{\sqrt{\lambda}}\right)\right\}\\
&=\left\{(t,x)\in[0,+\infty]\times[-\infty,+\infty]:\,\left(\frac{t}{\lambda},\frac{x}{\sqrt{\lambda}}\right)\in R_1\right\},
\end{align*}
where $R_1=\{(t,x)\in[0,+\infty]\times[-\infty,+\infty]:\,u_1(t,x)=v_{\mu_1}(x)\}$. Then, for every $(\lambda,x)\in]0,+\infty[\times[-\infty,+\infty]$,
$$
r_\lambda(x)=\inf\left\{t\in[0,+\infty]:\,\left(\frac{t}{\lambda},\frac{x}{\sqrt{\lambda}}\right)\in R_1\right\}=\lambda r_1\left(\frac{x}{\sqrt{\lambda}}\right)
$$
which completes the proof.  
\end{proof}
\begin{rem}
To prove that the function $(\lambda,x)\longmapsto r_\lambda(x)$ is self-similar, one may apply alternatively either the PDE characterisation of Root barrier provided by Gassiat-Oberhauser-Dos Reis (Theorem \ref{theo:RootFreeBounds} and Remark \ref{rem:RootBounds}) or, when the barrier function $r_1$ of $\mu_1$ is continuous, the integral equation for Root barrier due to Gassiat-Mijatovic-Oberhauser \cite{GMO15}.\\
1. Let $u_1$ be the unique viscosity solution of linear growth of the obstacle PDE
\begin{equation}\label{eq:RootPDE01}
\left\{
\begin{array}{ll}
\min\left(u-v_{\mu_{1}}, \partial_tu-\dfrac{1}{2}\partial_{xx} u\right)=0&\text{on }]0,+\infty[\times\R,\\
u(0,\cdot)=-|\cdot|&\text{on }\R,\\
u(+\infty,\cdot)=v_{\mu_{1}} &\text{on }\R,
\end{array}
\right.
\end{equation}
then, for every $\lambda>0$, the function $\widetilde{u}_\lambda$ defined by
$$
\widehat{u}_\lambda(t,x)=\sqrt{\lambda}u_1\left(\frac{t}{\lambda},\frac{x}{\sqrt{\lambda}}\right),
$$
is a viscosity solution of linear growth of (\ref{eq:RootPDE}), and one has
$$
\widehat{u}_\lambda(t,x)=u_\lambda(t,x),\quad\forall\,(t,x)\in [0,+\infty]\times[-\infty,+\infty].
$$
2. Suppose that $(\mu_\lambda,\lambda\geq0)$ is a family of atom-free probability measures such that, for every $\lambda>0$, $\mu_\lambda$ is the image of $\mu_1$ under $y\longmapsto\sqrt{\lambda}y$ and the regular barrier function $r_\lambda$ of the Root solution $\tau_\lambda$ for SEP$(\mu_\lambda)$ is continuous. This condition holds, for instance, when $\mu_1$ is  symmetric around $0$ and admits a compact support $[-\alpha,\alpha]$ and a bounded density which is nondecreasing on $[0,\alpha]$.  (see e.g. \cite[Proposition 1]{GMO15}).  Since $\mu_\lambda$ is also atom-free, its barrier function $r_\lambda$ solves the following Volterra integral equation:
\begin{align}
v_{\delta_0}(x)-v_{\mu_\lambda}(x)&=g\left(r_\lambda(x),x\right)-\nonumber\\
&\int_{\{y:\,r_\lambda(y)<r_\lambda(x)\}}g\left(r_\lambda(x)-r_\lambda(y),x-y\right)\mu_\lambda(dy)
\label{eq:RootVoltIntEql}\tag{EV$_\lambda$}
\end{align}
for every $x\in\R$, where $v_{\delta_0}=-|\cdot|$ and, for every $(t,z)\in\R_+\times\R$, 
$$
g(t,z)=\sqrt{\frac{2t}{\pi}}\exp\left(-\frac{z^2}{2t}\right)-|z|\text{Erfc}\left(\frac{|z|}{\sqrt{2t}}\right).
$$
Moreover, we know from corollary 2 in \citep{GMO15} that, for every $\lambda>0$, $r_\lambda$ is the unique continuous function that solves (\ref{eq:RootVoltIntEql}).
Observe that (\ref{eq:RootVoltIntEql}) is still valid if one replaces $x$ by $\sqrt{\lambda}x$. Then, as $\mu_\lambda$ is the image of $\mu_1$ under $y\longmapsto\sqrt{\lambda}y$, (\ref{eq:RootVoltIntEql}) rewrites
\begin{align*}
v_{\delta_0}(\sqrt{\lambda}x)-v_{\mu_\lambda}(\sqrt{\lambda}x)&=g\left(r_\lambda(\sqrt{\lambda}x),\sqrt{\lambda}x\right)-\nonumber\\
&\int_{\{y:\,r_\lambda(\sqrt{\lambda}y)<r_\lambda(\sqrt{\lambda}x)\}}g\left(r_\lambda(\sqrt{\lambda}x)-r_\lambda(\sqrt{\lambda}y),\sqrt{\lambda}x-\sqrt{\lambda}y\right)\mu_1(dy)
\label{eq:RootVoltIntEql1}\tag{EV$^{\prime}_\lambda$}
\end{align*}
But, since $(1/\sqrt{\lambda}) g(t,z)=g\left(t/\lambda,z/\sqrt{\lambda}\right)$, (\ref{eq:RootVoltIntEql1}) is equivalent to
\begin{align*}
v_{\delta_0}(x)-v_{\mu}(x)&=g\left(\widehat{r}_\lambda(x),x\right)-\\
&\int_{\{y:\,\widehat{r}_\lambda(y)<\widehat{r}(x)\}}g\left(\widehat{r}_\lambda(x)-\widehat{r}_\lambda(y),x-y\right)\mu_1(dy)\qquad\forall\,x\in\R,
\end{align*}
where, for every $z\in\R$, $\widehat{r}_\lambda(z)=(1/\lambda)r_\lambda(\sqrt{\lambda}z)$. Hence $\widehat{r}_\lambda$ is also a continuous function that solves (EV$_1$). It then follows from the uniqueness result  for (EV$_1$) that 
$$
r_1(x)=\widehat{r}_\lambda(x)=\frac{1}{\lambda}r_\lambda(\sqrt{\lambda}x),\quad\text{for all }x\in[-\infty,+\infty]
$$
which is equivalent to (\ref{eq:BarrierFunctEq}).
\end{rem}

\section{Root self-similar martingales}
In the next result we present a family of self-similar martingales. Precisely, we provide a necessary and sufficient condition so that the map $\lambda\longmapsto \tau_\lambda$ is non-decreasing a.s.. We also show that   several probability measures satisfy this condition.
\begin{theorem}\label{theo:Self-simRoot}
Let $\mu$ be an integrable and centered probability measure. Let $(B_v,v\geq0)$ be a standard Brownian motion started at $0$. For $\lambda\geq0$, let $\mu_{\lambda}$, $\tau_{\lambda}$ and $r_{\lambda}$ denote the image measure of $\mu$ under $x\longmapsto\sqrt{\lambda}x$,  the Root solution for SEP($\mu_{\lambda}$) and the regular barrier function of $\tau_{\lambda}$ respectively. 
\begin{enumerate}
\item[(i)]If the map $\lambda\longmapsto\tau_{\lambda}$ is a.s. non-decreasing, then 
$\left(B_{\tau_\lambda},\lambda\geq0\right)$ is a martingale satisfying Brownian scaling such that $B_{\tau_\lambda}$ has law $\mu_\lambda$ for every $\lambda\geq0$. Moreover, the process $\left((\tau_{\lambda},B_{\tau_{\lambda}}),\lambda\geq0\right)$ is Markovian and if, in addition, $\mu$ has no atom, then the martingale $\left(B_{\tau_\lambda},\lambda\geq0\right)$ is also Markovian.
\item[(ii)]The map $\lambda\longmapsto\tau_{\lambda}$ is a.s. non-decreasing if and only if 
\begin{equation}\label{eq:Self-simRoot}
x\longmapsto \dfrac{r_1(x)}{x^2}\text{ is non-decreasing on }]-\infty,0[\text{ and 
non-increasing on }]0,+\infty[.
\end{equation}
\end{enumerate} 
\end{theorem}
\begin{proof}\text{}
\item[(i)]Suppose that the map $\lambda\longmapsto\tau_\lambda$ is a.s. non-decreasing. Let $\left(\mathcal{F}_v,v\geq0\right)$ denote the natural filtration of $\left(B_v,v\geq0\right)$. Since, for every $\lambda\geq0$, $\left(B_{\tau_\lambda\wedge v},v\geq0\right)$ is uniformly integrable, then
$$
\E\left[B_{\tau_\eta}\left.\right|\mathcal{F}_{\tau_\lambda}\right]=B_{\tau_\lambda}\quad\text{for every }0\leq\lambda\leq\eta
$$
which means that $\left(B_{\tau_\lambda},\lambda\geq0\right)$ is a martingale. Moreover, by Brownian scaling, $$\left(W^{(c)}_v:=\frac{1}{c}B_{c^2v},v\geq0\right)$$ is still a one-dimensional Brownian motion for every $c>0$. As a consequence, we have
$$
\left(B_{\tau_\lambda},\lambda\geq0\right)\=\left(W^{(c)}_{\tau^{(c)}_\lambda},\lambda\geq0\right)=\left(\frac{1}{c}B_{c^2\tau^{(c)}_\lambda},\lambda\geq0\right),
$$
where
$$
\tau^{(c)}_\lambda=\inf\left\{v\geq0:\,v\geq r_\lambda\left(W_v^{(c)}\right)\right\}=\inf\left\{v\geq0:\,v\geq r_\lambda\left(\frac{B_{c^2v}}{c}\right)\right\}.
$$
Now, one may also deduce from Point (ii) of Theorem \ref{theo:RootPDEscale} that 
$$
r_\lambda\left(\frac{x}{c}\right)=\frac{1}{c^2}r_{c^2\lambda}(x),\quad\forall\,x\in[-\infty,+\infty].
$$
Indeed, one has
$$
r_\lambda\left(y\right)=\lambda r_1\left(\frac{y}{\sqrt{\lambda}}\right)=\frac{1}{c^2}r_{\lambda c^2}\left(c\,y\right),\quad\forall\,y\in[-\infty,+\infty].
$$
Hence,
$$
\tau^{(c)}_\lambda=\inf\left\{v\geq0:\,c^2v\geq r_{c^2\lambda}\left(B_{c^2v}\right)\right\}=\frac{1}{c^2}\tau_{c^2\lambda}
$$
and, as a consequence,
$$
\left(B_{\tau_\lambda},\lambda\geq0\right)\=\left(\frac{1}{c}B_{\tau_{c^2\lambda}},\lambda\geq0\right)
$$
which shows that $\left(B_{\tau_\lambda},\lambda\geq0\right)$ satisfies Brownian scaling.\\
Let $0<\lambda<\eta$ be fixed. Since $\tau_{\lambda}\leq\tau_{\eta}$ a.s., we have
\begin{align*}
\tau_{\eta}&=\inf\{t\geq\tau_{\lambda}:\,u_{\eta}(t,B_t)=v_{\mu_\eta}(B_t)\}\quad\text{a.s.}\\&=\tau_{\lambda}+\inf\left\{s\geq0:\,u_{\eta}\left(s+\tau_{\lambda},B_{s+\tau_{\lambda}}\right)= v_{\mu_\eta}(B_{s+\tau_{\lambda}})\right\}
\\&=
\tau_{\lambda}+\theta_{\eta}\left(\tau_{\lambda},B^{\tau_{\lambda},B_{\tau_\lambda}}_{\cdot+\tau_{\lambda}}:=\left(B^{\tau_{\lambda},B_{\tau_{\lambda}}}_{v+\tau_{\lambda}},v\geq0\right)\right),
\end{align*}
where
$$
B^{s,x}_t=x+B_{t}-B_{s},\quad\text{for all }(s,x)\in\R_+\times\R,\text{ }t\geq s,
$$
and
$$
\theta_{\eta}\left(\tau_{\lambda},B^{\tau_{\lambda},B_{\tau_\lambda}}_{\cdot+\tau_{\lambda}}\right):=\inf\left\{s\geq0:\,u_{\eta}\left(s+\tau_{\lambda},B_{\tau_\lambda}+B^{\tau_\lambda,0}_{s+\tau_\lambda}\right)=v_{\mu_\eta}\left(B_{\tau_\lambda}+B^{\tau_\lambda,0}_{s+\tau_\lambda}\right)\right\}.
$$
It follows from the strong Markov property of Brownian motion that $B^{\tau_\lambda,0}_{\cdot+\tau_\lambda}=\left(B^{\tau_\lambda,0}_{v+\tau_\lambda},v\geq0\right)$
is a Brownian motion independent of $\mathcal{F}_{\tau_\lambda}$. Then, for every bounded measurable function $\phi$,
$$
\E\left[\phi(\tau_{\eta},B_{\tau_{\eta}})\left\vert\mathcal{F}_{\tau_{\lambda}}\right.\right]=\psi\left(\tau_{\lambda},B_{\tau_{\lambda}}\right),
$$
where
$$
\psi(t,x)=\E\left[\phi\left(t+\theta_{\eta}\left(t,B^{t,x}_{\cdot+t}\right),B^{t,x}_{t+\theta_{\eta}\left(t,B^{t,x}_{\cdot+t}\right)}\right)\right]
$$
which shows that $\left((\tau_\lambda,B_{\tau_{\lambda}}),\lambda\geq0\right)$ is a (non-homogeneous) Markov process.
If $\mu$ has no atom, then one deduces from Lemma 1 in \citep{GMO15} that $r(B_{\tau_{\lambda}})=\tau_{\lambda}$ a.s. and, as a consequence, that $\left(B_{\tau_\lambda},\lambda\geq0\right)$ is Markovian.  

\item[(ii)]Let $R_\lambda$ denote the Root barrier given by the function $r_\lambda$. We first note that the map $\lambda\longmapsto \tau_\lambda$ is a.s. non-decreasing if and only if the family $(R_\lambda,\lambda\geq0)$ is non-increasing in the sense of set inclusion, which means that $\lambda\longmapsto r_\lambda(x)$ is non-decreasing for every $x\in\R$. Hence it remains to show that $\lambda\longmapsto r_\lambda(x)$ for every $x\in\R$ if and only if $x\longmapsto r_1(x)/x^2$ is non-decreasing on $]-\infty,0[$ and non-increasing on $]0,+\infty[$. Suppose first that $\lambda\longmapsto r_\lambda(x)=\lambda r_1\left(x/\sqrt{\lambda}\right)$ is non-decreasing for every $x\in\R$. Let $x_1\leq x_2$ belong to $]-\infty,0[$. Since $1\geq\lambda=x_2^2/x_1^2$, one has
$$
r_1(x_2)\geq\lambda r_1\left(\frac{x_2}{\sqrt{\lambda}}\right)=x_2^2\frac{r_1(x_1)}{x_1^2}
$$
which shows that $x\longmapsto r_1(x)/x^2$ is non-decreasing on $]-\infty,0[$. Similarly, one may show that $x\longmapsto r_1(x)/x^2$ is non-increasing on $]0,+\infty[$. Conversely, suppose that $x\longmapsto r_1(x)/x^2$ is non-decreasing on $]-\infty,0[$ and non-increasing on $]0,+\infty[$. As $r(0)$ is positive, $\lambda\longmapsto \lambda r_1(0)$ is non-decreasing. Moreover, for every $x\in\R\setminus\{0\}$, one may observe that
$$
\lambda r_1\left(\frac{x}{\sqrt{\lambda}}\right)=x^2\frac{r_1\left(x/\sqrt{\lambda}\right)}{x^2/\lambda}
$$
which shows that $\lambda\longmapsto\lambda r_1\left(x/\sqrt{\lambda}\right)$ is still non-decreasing when $x\in\R\setminus\{0\}$. 
\end{proof}
\begin{rem}
Condition (\ref{eq:Self-simRoot}) is weaker than
\begin{equation}\label{eq:Self-simRoot1}
x\longmapsto r_1(x)\text{ is non-decreasing on }]-\infty,0[\text{ and 
non-increasing on }]0,+\infty[.
\end{equation}
Indeed, if (\ref{eq:Self-simRoot1}) holds, then, as $x\longmapsto 1/x^2$ is positive, non-decreasing on $]-\infty,0[$, non-increasing on $]0,+\infty[$ and as $r_1$ is non-negative, one deduces that $x\longmapsto r_1(x)/x^2$ is non-decreasing on $]-\infty,0[$, and non-increasing on $]0,+\infty[$.
Many regular barrier functions given in the literature satisfy Condition (\ref{eq:Self-simRoot1}) (see e.g. \citep[Example 5.2]{Ho11} and \citep[Section 2]{CW13}).   
\end{rem}
We now present some examples to which Theorem \ref{theo:Self-simRoot} applies. We mention that, given a probability measure $\mu_1$, it is hard to compute the regular barrier function $r_1$ that gives the Root Solution to SEP($\mu_1$). There are only few cases where this can be done explicitely. Numerical methods have been provided to compute Root barriers with great  precision (see e.g. \cite[Section 3]{GMO15} and  \cite[Section 4]{GODR17}).
\begin{exa}\label{exa:RootBarFunct}
\item[1.]If $\mu_1$ is a zero-mean normal distribution, then it is not difficult to see that $r_1$ is constant. Indeed, the Root solution $\tau_1$ to SEP($\mu_1$) equals to the square-mean of $\mu_1$. In this case $(B_{\tau_\lambda},\lambda\geq0)$ is still a Brownian motion which may be non-standard.
\item[2.]If $\mu_1$ is of the form
$$
\mu_1=\frac{b}{a+b}\delta_{a}-\frac{a}{a+b}\delta_b,
$$
where $a,b$ are real numbers such that $a<0<b$, then the corresponding Root barrier function is (see e.g. \cite[Section 2]{CW13})
$$
r_1(x)=\left\{
\begin{array}{ll}
\infty &\text{if }x\in]a,b[\\
0 &\text{if }x\notin ]a,b[.
\end{array}
\right.
$$
\item[3.]Suppose that $\mu_1$ has the form 
$$
\mu_1=p\delta_{-a}+(1-2p)\delta_0+p\delta_{a}\quad(1\leq p\leq 1/2,\,a>0).
$$ 
Then the corresponding regular barrier function is 
$$
r_1(x)=\left\{
\begin{array}{cl}
0&\text{if }x\in]-\infty,-a]\cup[a,+\infty[\\
\infty&\text{if }x\in]-a,0[\cup]0,a[\\
t_{0}(a,p)&\text{if }x=0,
\end{array}
\right.
$$
where $t_{0}(a,p)$ is a nonnegative real number (see e.g. \cite[Point 3 of Example 5.2]{Ho11}). Observe that $r_1$ is non-decreasing on $]-\infty,0[$ and non-increasing on $]0,+\infty[$ . 
\item[4.]Let $\mu_1$ be an integrable probability measure on $(\R,\mathcal{B}(\R))$ which satisfies the Assumption 1 in \citep{GMO15}, that is $\mu_1$ has mean zero, and the Root barrier solving SEP($\mu$) is given by a function $r_1$ which is symmetric around $0$, continuous, and non-increasing on $[0,+\infty]$.  There are several probability measures satisfying the above assumption. Indeed, as proved in \cite[Proposition 1]{GMO15}, symmetric probability measures around $0$ with compact support $[-\alpha,\alpha]$ and bounded non-decreasing density on $[0,\alpha]$ fulfill Assumption 1 in \cite{GMO15}.
\end{exa}
We mention that there are also probability measures to which  Theorem \ref{theo:Self-simRoot} does not apply. For instance, if $\mu_1$ is the canonical measure on the middle Cantor set, then, by Root's result, the resulting barrier function $r_1$ must be finite only on the Cantor set.
 
\section{Numerics: Monotone self-similar Root barriers}
The aim of this paragraph is to give some pictorial representations of the monotonicity property of certain self-similar regular Root barriers. We use an explicit finite differences scheme adapted from the scheme implemented in \cite[Section 4]{GODR17} to simulate Root barriers. The Barles-Souganidis method \cite{BS90,BS91} gives a convergence result of this scheme.
\\
Fix an integrable mean-zero probability measure $\mu_1$ whose support is denoted by supp$(\mu_1)$. Fix also $\lambda\in]0,1[$ and $T>0$. Let $\mu_\lambda$ be the image of $\mu_1$ under $y\longmapsto\sqrt{\lambda}y$. We distinguish two cases.
%{\bf Case where supp$(\mu_1)$ is bounded.} 
Choose $a<0<b$ such that supp$(\nu_1)\subset[a,b]$, and consider the following time-space mesh of points in $\mathcal{O}_T:=[0,T]\times[a,b]$
$$
\mathcal{G}_h:=\left\{t_n:\,t_n=n\Delta t,\,n=0,1,\cdots,N_T\right\}\times\left\{x_j:\,x_j=a+j\Delta x,\,j=0,1,\cdots,N_x\right\},
$$
where $h=(\Delta t,\Delta x)=(T/N_T,(b-a)/N_x)$ with $N_T,\,N_x\in\N$ large enough. As in \cite[Section 4]{GODR17}, we denote by $\mathcal{B}(\mathcal{O}_T,\R)$ the set of bounded function from $\mathcal{O}_T$ to $\R$ and by $\mathcal{B}\mathcal{U}\mathcal{C}(\mathcal{O}_T,\R)$ the subset of $\mathcal{B}(\mathcal{O}_T,\R)$ consists of bounded uniformly continuous functions. Since supp$(\mu_1)$ is bounded, the unique viscosity solution $u_1$, resp. $u_\lambda$ of (\ref{eq:RootPDE01}), resp. (\ref{eq:RootPDE}) belongs to $\mathcal{B}\mathcal{U}\mathcal{C}(\mathcal{O}_T,\R)$. 
Let $u_{0,\alpha}^h:\,\mathcal{G}_h\to\R$  (with  $\alpha\in\{\lambda,1\}$) denote the discrete approximation of $u_{\alpha}$. The values of $u_{0,\alpha}^h$ are obtained by solving the system
$$
\left\{
\begin{array}{ll}
u^h_{0,\alpha}(0,x_j)=-|x_j|,&\text{for }j=0,1,\cdots,N_x,\\
u^h_{0,\alpha}(t_{n+1},x_j)=\max\left\{v_{\mu_{\alpha}}(x_j),S_0^h[u_{0,\alpha}^h](t_n,x_j)\right\},&\text{for }(n,j)\in[0,N_T-1]\times[0,N_x].
\end{array}
\right.
$$
with $S_0^h\left[u^h_{0,\alpha}\right]$ defined as
\begin{align*}
&S_0^h\left[u^h_{0,\alpha}\right](t_n,x_j):=\\
&\left\{
\begin{array}{cl}
-|x_j|&\text{if }n=0,\\
-|x_j|&\text{if }j\in\{0,N_x\},\\
u_{0,\alpha}^h(t_n,x_j)+\frac{\Delta t}{2(\Delta x)^2}\left(u_{0,\alpha}^h(t_n,x_{j+1})-2u_{0,\alpha}^h(t_n,x_j)+u_{0,\alpha}^h(t_n,x_{j-1})\right)&\text{otherwise},
\end{array}
\right.
\end{align*}
where we suppose that the usual CFL condition: $\Delta t<(\Delta x)^2$ holds.
Let $u_{\alpha}^h:\,\mathcal{O}^h_T\to\R$ (with $\mathcal{O}_T^h:=[0,T+\Delta t[\times[a-\Delta x/2,b+\Delta x/2]$) be the function  given by $u_{\alpha}^h(t,x)=u_{0,\alpha}^h(t_n,x_j)$ when $(t,x)\in[t_n,t_{n+1}[\times[x_j-\Delta x/2,x_j+\Delta x/2[$ for some $n\in\{0,1,\cdots,N_T\}$ and $j\in\{0,1,\cdots,N_x\}$. Observe that the restriction of $u_{\alpha}^h$ to $\mathcal{O}_T$, which is also denoted by $u_{\alpha}^h$, belongs to $\mathcal{B}(\mathcal{O}_T,\R)$. Since $\mu_{\alpha}$ has bounded support, $u^h_{\alpha}\in\mathcal{B}([0,T]\times\R,\R)$. By Proposition 1 in \cite{GODR17}, 
\begin{equation}\label{eq:ApproxiConv}
\lim\limits_{h\to(0,0)}\sup\limits_{[0,T]\times\R}\left|u^h_{\alpha}-u_{\alpha}\right|=0.
\end{equation}
Moreover, Proposition 2 in \cite{GODR17} provides the rate of the convergence result (\ref{eq:ApproxiConv}) (see \cite[Section 3]{Ja03} for more details). Hence, for sufficiently small $h$, the subset $\mathcal{R}^h_{\alpha}$ of $\mathcal{O}_T$, defined as
\begin{equation}
\mathcal{R}^h_{\alpha}:=\left\{(t,x)\in\mathcal{O}_T:\,u_{\alpha}^h(t,x)=v_{\mu_{\alpha}}(x)\right\},
\end{equation}
nearly coincides with the regular Root barrier $\mathcal{R}_{\alpha}$ on the rectangle $\mathcal{O}_T$.  

\begin{figure}[h!]
\begin{center}
\includegraphics[width=4.00cm]{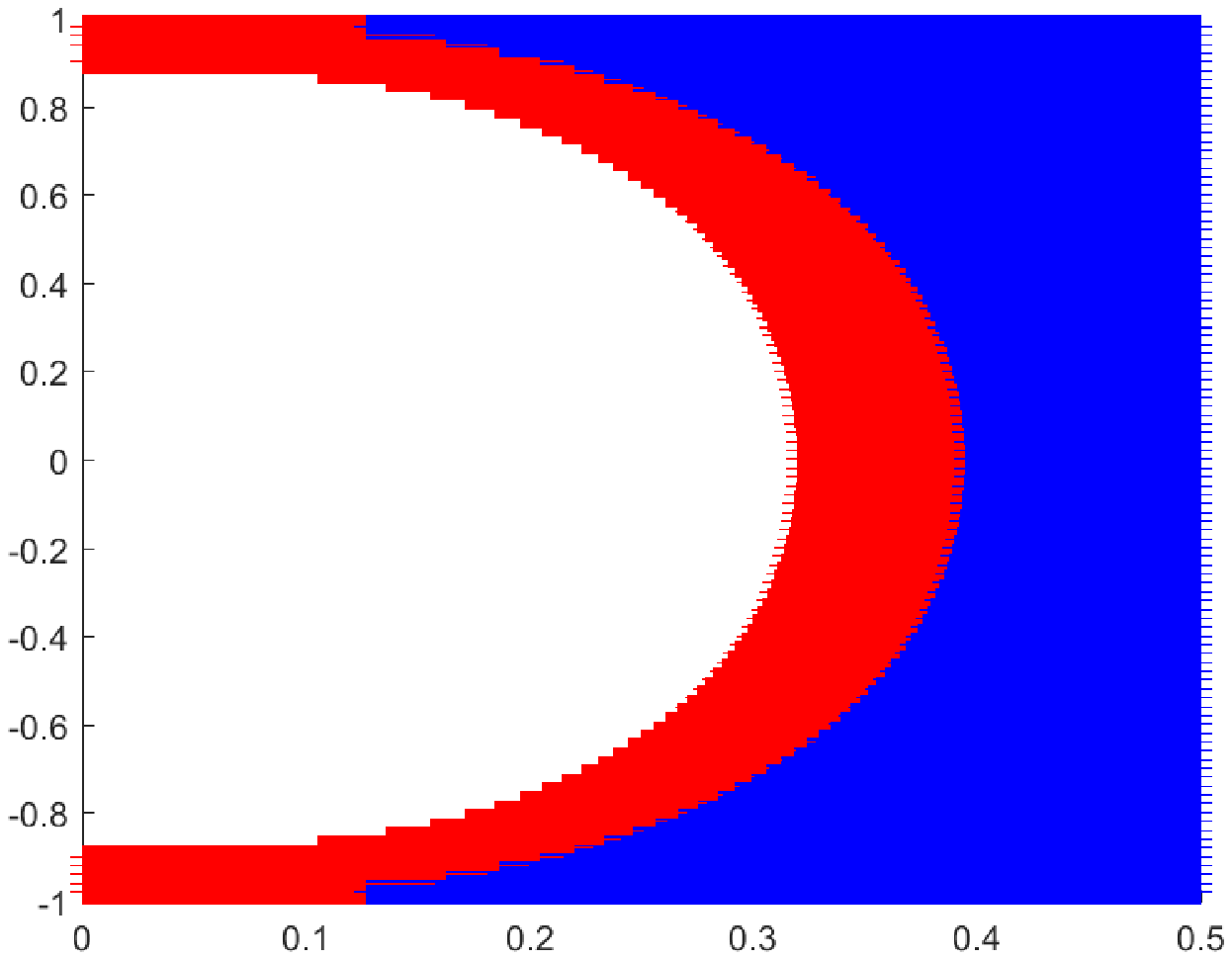}
\includegraphics[width=4.00cm]{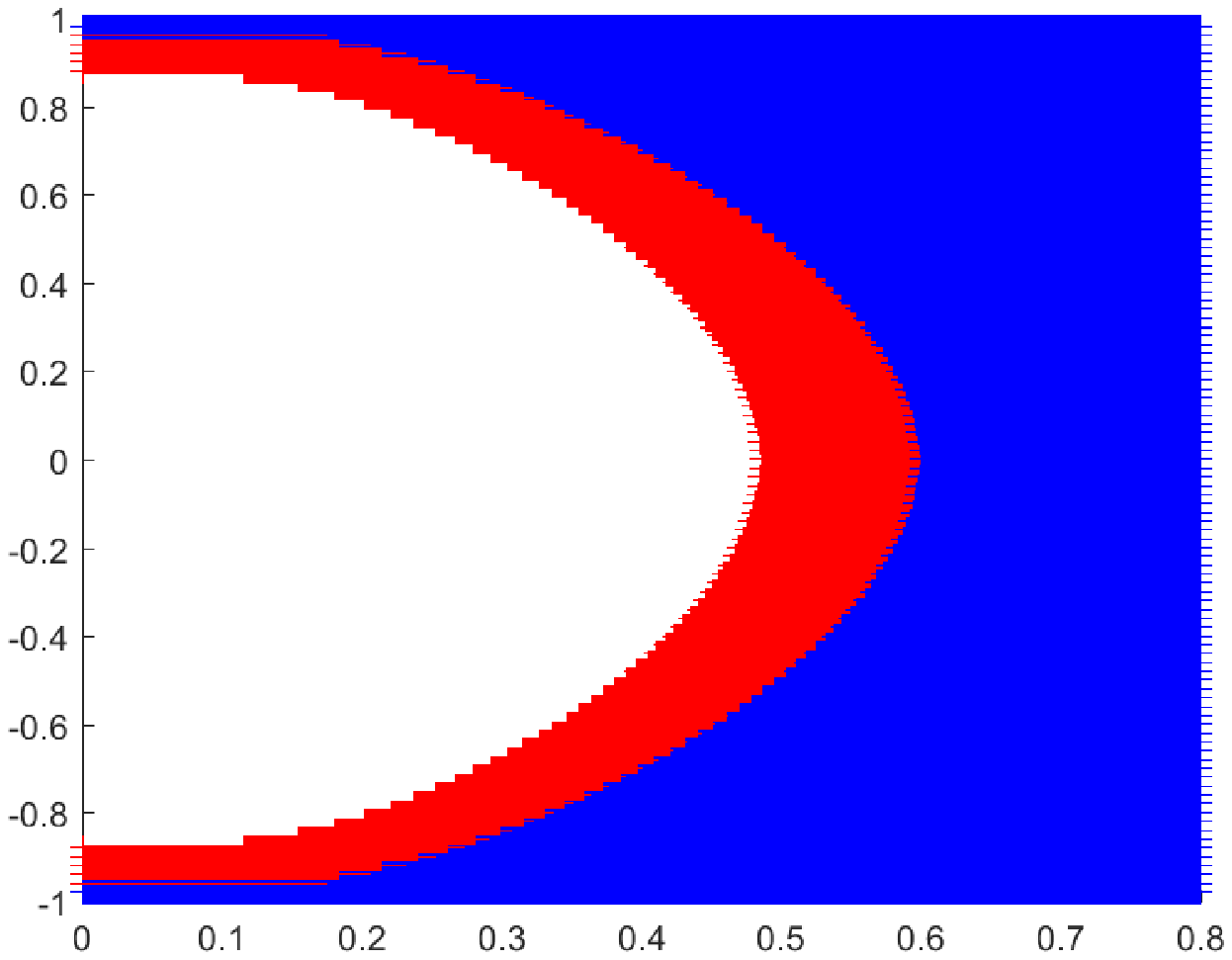}
%,height=5.500cm
\includegraphics[width=4.00cm]{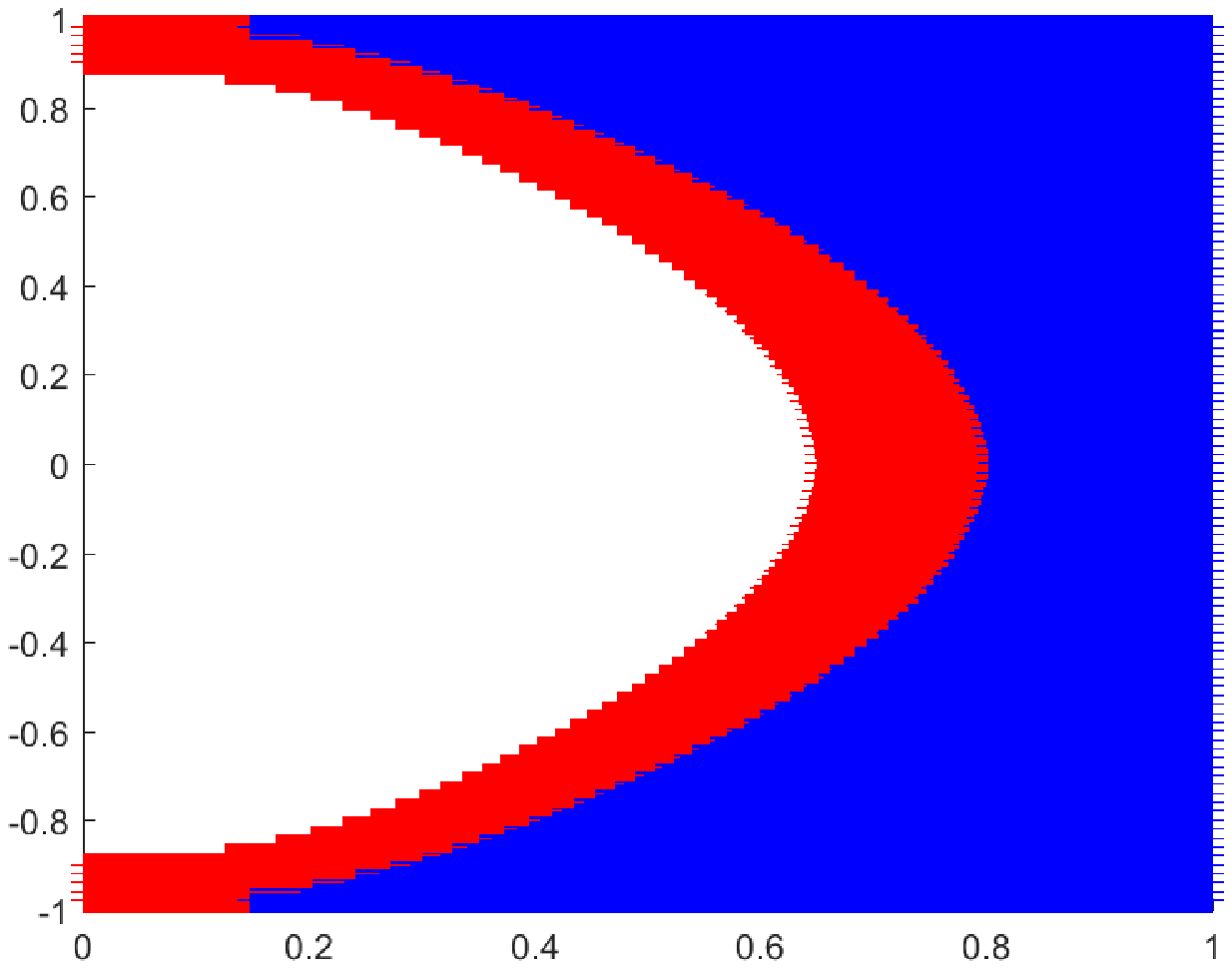}
%drivencavity.eps}
\end{center}
\caption{$\mu_1(dx)=\mathbf{1}_{[-1,1]}(x)dx$ (left plot), $\mu_1(dx)=0.75\sqrt{|x|}\mathbf{1}_{[-1,1]}dx$ (middle plot), $\mu_1(dx)=|x|\mathbf{1}_{[-1,1]}(x)dx$ (right plot), $\lambda=0.81$, $h=0.02$, $N_T=20000$.}\label{Fig:figure1}
\end{figure}

\begin{figure}[h!]
\begin{center}
\includegraphics[width=4.00cm]{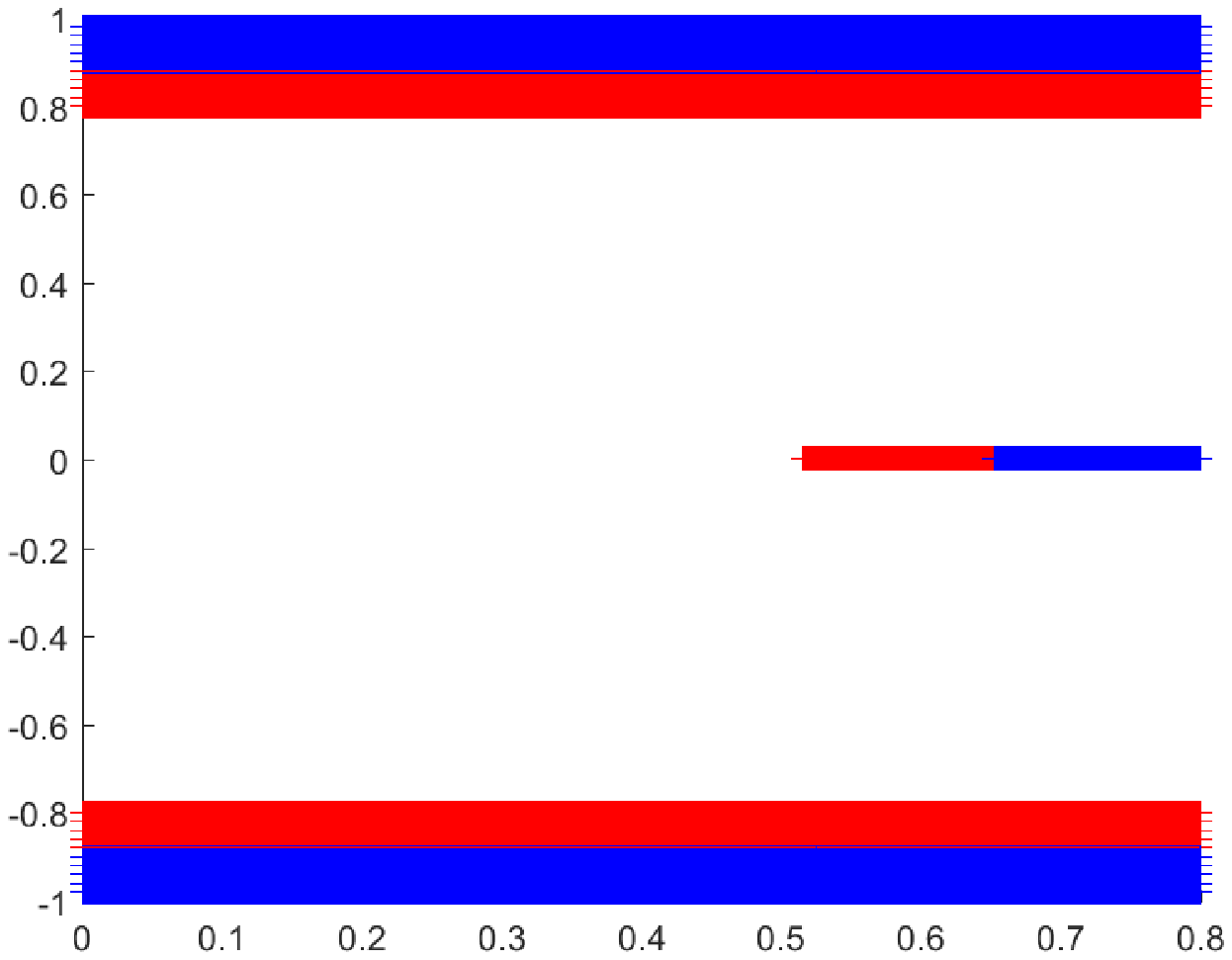}
\includegraphics[width=4.00cm]{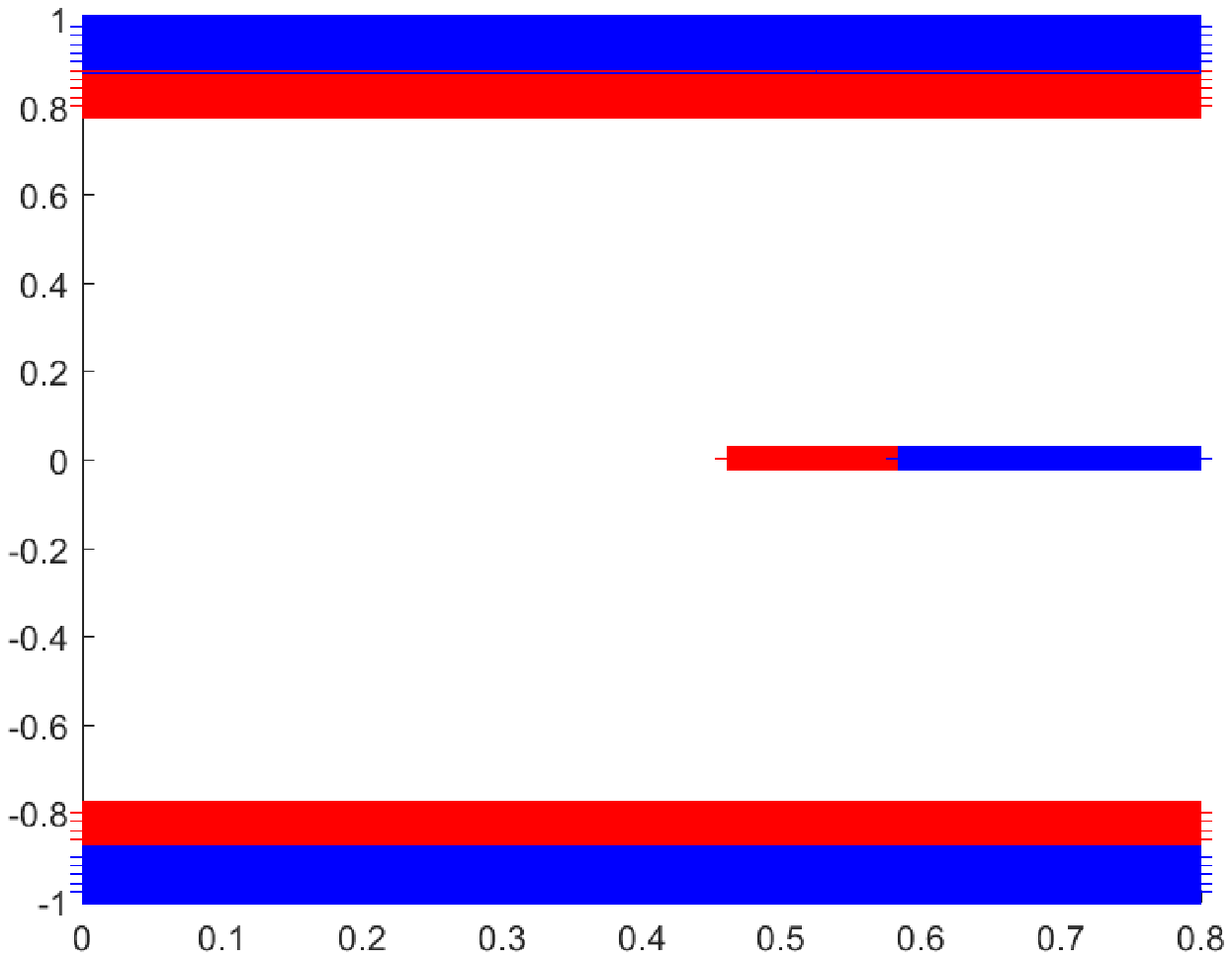}
\includegraphics[width=4.00cm]{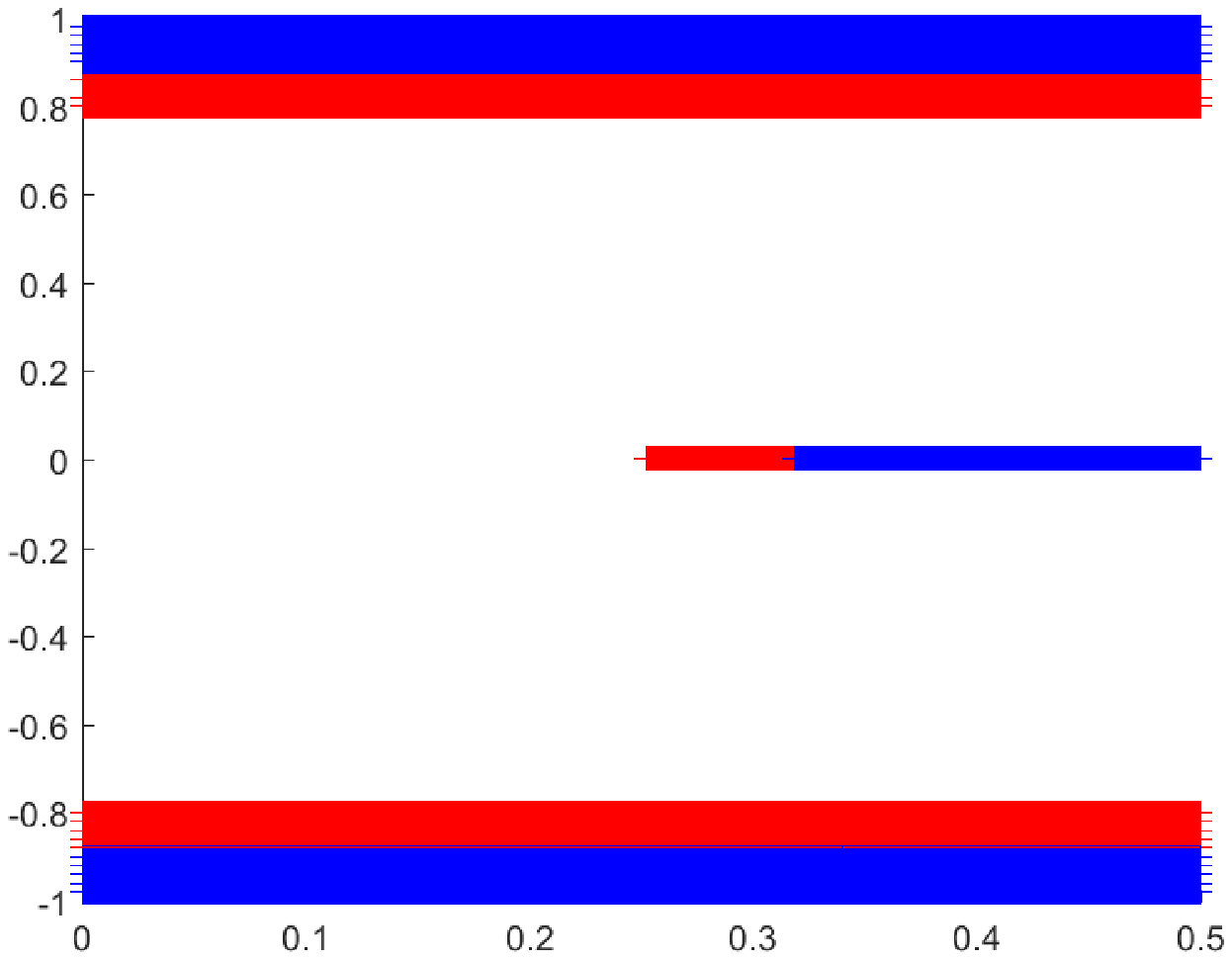}
\end{center}
\caption{$\mu_1=(7/20)\delta_{-0.9}+(3/10)\delta_0+(7/20)\delta_{0.9}$ (left plot), $\mu_1=(1/3)\delta_{-0.9}+(1/3)\delta_0+(1/3)\delta_{0.9}$ (middle plot), $\mu_1=(1/4)\delta_{-0.9}+(1/2)\delta_0+(1/4)\delta_{0.9}$ (right plot), $\lambda=\frac{0.64}{0.81}$, $h=0.02$, $N_T=20000$.}\label{Fig:figure2}
\end{figure}

The figures below show that, for $h$ sufficiently small, $\mathcal{R}^h_{\lambda}$ includes $\mathcal{R}^h_1$  when $\mu_1$ satisfies Condition (\ref{eq:Self-simRoot}) in Theorem \ref{theo:Self-simRoot}. The barriers $\mathcal{R}^h_1$ and $\mathcal{R}^h_{\lambda}$ are plotted in blue and red respectively. 
Figure \ref{Fig:figure1} illustrate the fourth point in Example \ref{exa:RootBarFunct}. Precisely, $\mu_1$ is a symmetric probability measure around $0$ with compact support $[-1,1]$ and a bounded non-decreasing density on $[0,1]$. Figure \ref{Fig:figure2} is a pictorial representation of the third point of Example \ref{exa:RootBarFunct}. Indeed, $\mu_1$ has the form
$$
\mu_1=p\delta_{-a}+(1-2p)\delta_0+p\delta_{a},
$$
where $0<p<1/2$ and $a>0$. 

%\bibliographystyle{apalike}
%\bibliography{Root-Embedding}

\newpage
\begin{thebibliography}{}

\bibitem[Albin, 2008]{Al08}
Albin, J. M.~P. (2008).
\newblock A continuous non-brownian motion martingale with brownian motion
  marginal distributions.
\newblock {\em Statist. Probab. Letters}, 78(6):682--686.

\bibitem[Baker et~al., 2011]{BDY11}
Baker, D., Donati-Martin, C., and Yor, M. (2011).
\newblock A sequence of albin type continuous martingales with brownian
  marginals and scaling.
\newblock In {\em S{\'e}minaire de probabilit{\'e}s XLIII}, pages 441--449.
  Springer.

\bibitem[{Barles} and {Souganidis}, 1990]{BS90}
{Barles}, G. and {Souganidis}, P.~E. (1990).
\newblock Convergence of approximation schemes for fully nonlinear second order
  equations.
\newblock In {\em 29th IEEE Conference on Decision and Control}, pages
  2347--2349 vol.4.

\bibitem[Barles and Souganidis, 1991]{BS91}
Barles, G. and Souganidis, P.~E. (1991).
\newblock Convergence of approximation schemes for fully nonlinear second order
  equations.
\newblock {\em Asymptotic analysis}, 4(3):271--283.

\bibitem[{Beiglb{\"o}ck} et~al., 2017]{BCH17}
{Beiglb{\"o}ck}, M., Cox, A. M.~G., and Huesmann, M. (2017).
\newblock Optimal transport and skorokhod embedding.
\newblock {\em Inventiones mathematicae}, 208(2):327--400.

\bibitem[Bogso, 2015]{Bo15}
Bogso, A.~M. (2015).
\newblock Mrl order, log-concavity and an application to peacocks.
\newblock {\em Stochastic Processes and their Applications}, 125(4):1282--1306.

\bibitem[Cox et~al., 2017]{COT17}
Cox, A.~M., Ob{\l}{\'o}j, J., and Touzi, N. (2017).
\newblock The root solution to the multi-marginal embedding problem: an optimal
  stopping and time-reversal approach.
\newblock {\em arXiv preprint arXiv:1505.03169v2}.

\bibitem[Cox et~al., 2013]{CW13}
Cox, A.~M., Wang, J., et~al. (2013).
\newblock Root’s barrier: Construction, optimality and applications to
  variance options.
\newblock {\em The Annals of Applied Probability}, 23(3):859--894.

\bibitem[Dupire, 2005]{Du05}
Dupire, B. (2005).
\newblock Arbitrage bounds for volatility derivatives as free boundary problem.
\newblock {\em Presentation at PDE and Mathematical Finance, KTH, Stockholm}.

\bibitem[Fan et~al., 2015]{FHK15}
Fan, J.~Y., Hamza, K., Klebaner, F., et~al. (2015).
\newblock Mimicking self-similar processes.
\newblock {\em Bernoulli}, 21(3):1341--1360.

\bibitem[Gassiat et~al., 2015]{GMO15}
Gassiat, P., Mijatovi{\'c}, A., Oberhauser, H., et~al. (2015).
\newblock An integral equation for root’s barrier and the generation of
  brownian increments.
\newblock {\em The Annals of Applied Probability}, 25(4):2039--2065.

\bibitem[Gassiat et~al., 2017]{GODR17}
Gassiat, P., Oberhauser, H., and dos Reis, G. (2017).
\newblock Root’s barrier, viscosity solutions of obstacle problems and
  reflected fbsdes.
\newblock {\em Stochastic Process. Appl.}, 125(12):4601--4631.

\bibitem[Hamza and Klebaner, 2007]{HK07}
Hamza, K. and Klebaner, F.~C. (2007).
\newblock A family of non-gaussian martingales with gaussian marginals.
\newblock {\em International Journal of Stochastic Analysis}, 2007.

\bibitem[Henry-Labordere et~al., 2016]{HLTT16}
Henry-Labordere, P., Tan, X., and Touzi, N. (2016).
\newblock An explicit martingale version of the one-dimensional brenier’s
  theorem with full marginals constraint.
\newblock {\em Stochastic Process. Appl.}, 126(9):2800--2834.

\bibitem[Hirsch et~al., 2011]{HPRY11}
Hirsch, F., Profeta, C., Roynette, B., and Yor, M. (2011).
\newblock Constructing self-similar martingales via two skorokhod embeddings.
\newblock In {\em S{\'e}minaire de probabilit{\'e}s XLIII}, pages 451--503.
  Springer.

\bibitem[Hobson, 2011]{Ho11}
Hobson, D. (2011).
\newblock The skorokhod embedding problem and model-independent bounds for
  option prices.
\newblock In {\em Paris-Princeton Lectures on Mathematical Finance 2010}, pages
  267--318. Springer.

\bibitem[Hobson, 2013]{Ho13}
Hobson, D.~G. (2013).
\newblock Fake exponential brownian motion.
\newblock {\em Statist. Probab. Letters}, 83(10):2386--2390.

\bibitem[Jakobsen, 2003]{Ja03}
Jakobsen, E.~R. (2003).
\newblock On the rate of convergence of approximation schemes for bellman
  equations associated with optimal stopping time problems.
\newblock {\em Mathematical Models and Methods in Applied Sciences},
  13(05):613--644.

\bibitem[{Jourdain} and {Zhou}, 2016]{JZ16}
{Jourdain}, B. and {Zhou}, A. (2016).
\newblock {Existence of a calibrated regime switching local volatility model
  and new fake Brownian motions}.
\newblock {\em ArXiv e-prints}.

\bibitem[K{\"a}llblad et~al., 2017]{KTT17}
K{\"a}llblad, S., Tan, X., and Touzi, N. (2017).
\newblock Optimal skorokhod embedding given full marginals and az{\'e}ma--yor
  peacocks.
\newblock {\em The Annals of Applied Probability}, 27(2):686--719.

\bibitem[Kiefer, 1972]{Kie72}
Kiefer, J. (1972).
\newblock Skorohod embedding of multivariate rv's, and the sample df.
\newblock {\em Probability Theory and Related Fields}, 24(1):1--35.

\bibitem[Loynes, 1970]{Lo70}
Loynes, R.~M. (1970).
\newblock Stopping times on brownian motion: Some properties of root's
  construction.
\newblock {\em Probab. Theory Related Fields}, 16(3):211--218.

\bibitem[Madan and Yor, 2002]{MY02}
Madan, D.~B. and Yor, M. (2002).
\newblock Making markov martingales meet marginals: with explicit
  constructions.
\newblock {\em Bernoulli}, 8(4):509--536.

\bibitem[Oleszkiewicz, 2008]{Ol08}
Oleszkiewicz, K. (2008).
\newblock On fake brownian motions.
\newblock {\em Statist. Probab. Letters}, 78(11):1251--1254.

\bibitem[Richard et~al., 2018]{RTT18}
Richard, A., Tan, X., and Touzi, N. (2018).
\newblock On the root solution to the skorokhod embedding problem given full
  marginals.
\newblock {\em arXiv preprint arXiv:1810.10048}.

\bibitem[Root, 1969]{Ro69}
Root, D.~H. (1969).
\newblock The existence of certain stopping times on brownian motion.
\newblock {\em Ann. Math. Statist.}, 40(2):715--718.

\bibitem[Rost, 1976]{Rst76}
Rost, H. (1976).
\newblock Skorokhod stopping times of minimal variance.
\newblock In {\em S{\'e}minaire de Probabilit{\'e}s X Universit{\'e} de
  Strasbourg}, pages 194--208. Springer.

\bibitem[Skorokhod, 2014]{Sk14}
Skorokhod, A.~V. (2014).
\newblock {\em Studies in the theory of random processes}.
\newblock Courier Corporation.

\end{thebibliography}

\end{document}